\documentclass[11pt]{article}
\usepackage{latexsym,amssymb,amsmath}
\def\demo{\noindent{\bf Proof. }}
\def\QED{\hfill$\Box$}

\newtheorem{Theorem}{Theorem}[section]
\newtheorem{Lemma}[Theorem]{Lemma}
\newtheorem{Corollary}[Theorem]{Corollary}
\newtheorem{Proposition}[Theorem]{Proposition}
\newtheorem{Remark}[Theorem]{Remark}
\newtheorem{Example}[Theorem]{Example}

\newtheorem{Definition}[Theorem]{Definition}

\begin{document}

\topmargin3mm
\hoffset=-1cm
\voffset=-1.5cm
\begin{center}

{\large\bf Rees algebras and polyhedral cones of ideals 
of \\ vertex
covers of  perfect graphs}\\

\vspace{6mm}


\addtocounter{footnote}{1}

\medskip

Rafael
H. Villarreal\footnote{Partially supported
by CONACyT grant 49251-F and SNI, M\'exico.}\\

{\small Departamento de
Matem\'aticas}\vspace{-1mm}\\
{\small Centro de Investigaci\'on y de Estudios
Avanzados del
IPN}\vspace{-1mm}\\
{\small Apartado Postal
14--740}\vspace{-1mm}\\
{\small 07000 Mexico City, D.F.}\vspace{-1mm}\\
{\small
e-mail: {\tt
vila@math.cinvestav.mx}}\vspace{4mm}

\end{center}

\date{}

\begin{abstract}
\noindent Let $G$ be a perfect graph and let $J$ be its ideal of
vertex covers. We show that the Rees algebra of $J$ is
normal and that this algebra is Gorenstein if $G$ is unmixed. 
Then we
give a description--in terms of cliques--of the 
symbolic Rees algebra and the Simis cone of the edge ideal of $G$. 
\end{abstract} 

\section{Introduction} 

Let $R=K[x_1,\ldots,x_n]$ be a polynomial ring 
over a field $K$ and let $I$ be an ideal 
of $R$ of height $g\geq 2$ minimally generated by a finite set
$F=\{x^{v_1},\ldots,x^{v_q}\}$ of
square-free monomials  
 of degree at least two. As usual we
use  
$x^a$ as an abbreviation for $x_1^{a_1} \cdots x_n^{a_n}$, 
where $a=(a_1,\ldots,a_n)\in \mathbb{N}^n$. A {\it clutter\/} with  
vertex set  $X$ is a family of subsets of $X$,
called edges, none of which is included in another. The set of
vertices and edges of $\cal C$ are denoted by $V({\cal C})$ and
$E({\cal C})$ respectively.    
We can associate to the 
ideal $I$ a {\it clutter\/} $\cal C$ by taking the set 
of indeterminates $X=\{x_1,\ldots,x_n\}$ as vertex set and 
$E=\{S_1,\ldots,S_q\}$ as edge set, where $S_k$ is the support of 
$x^{v_k}$, i.e., $S_k$ is the set of variables that occur in
$x^{v_k}$.  For this
reason  $I$ is called the {\it edge ideal\/} of $\cal C$. To stress 
the 
relationship  between $I$ and $\cal C$ we will use the notation
$I=I({\cal C})$.  The $n\times q$ matrix with column vectors
$v_1,\ldots,v_q$ will be denoted by $A$, it is called the {\it
incidence 
matrix\/} of $\cal C$. It is usual to call $v_i$ the {\it 
incidence vector\/} or {\it characteristic vector\/} of $S_i$. 

The {\it blowup algebras\/} studied here are: 
(a) the {\it Rees algebra\/}
$$
R[It]=R\oplus It\oplus\cdots\oplus I^{i}t^i\oplus\cdots
\subset R[t],
$$
where $t$ is a new variable, and (b) 
the {\it symbolic Rees algebra}
$$
R_s(I)=R\oplus I^{(1)}t\oplus\cdots\oplus
I^{(i)}t^i\oplus\cdots\subset R[t], 
$$
where $I^{(i)}$ is the 
$i${\it th\/} symbolic power of $I$. 

The {\it Rees cone\/} of $I$, denoted by $\mathbb{R}_+(I)$, is the 
polyhedral  
cone consisting  
of the non-negative linear combinations of the set 
$$
{\cal
A}'=\{e_1,\ldots,e_n,(v_1,1),\ldots,(v_q,1)\}\subset\mathbb{R}^{n+1}, 
$$

where $e_i$ is the $i${\it th} unit vector. It is well documented
\cite{normali,reesclu,clutters} that 
Rees cones are an  
effective device to study algebraic and combinatorial properties of
blowup algebras of square-free monomial ideals and clutters. They 
will play an important role here (Lemma~\ref{mundial09-jun-06}). The 
normalization of $R[It]$ can be expressed in terms of Rees cones as
we now explain. Let $\mathbb{N}{\cal A}'$ be the subsemigroup of 
$\mathbb{N}^{n+1}$ generated by ${\cal A}'$, consisting of the 
linear combinations of ${\cal A}'$ with non-negative integer
coefficients. The Rees algebra of $I$ can be written as
\begin{equation}\label{may6-06-2}
R[It]=
K[\{x^at^b\vert\, (a,b)\in\mathbb{N}{\cal A}'\}].
\end{equation}
According to \cite[Theorem~7.2.28]{monalg} the  
{\it integral closure\/} of $R[It]$ 
in its field of fractions can be expressed as
\begin{equation}\label{may6-06}
\overline{R[It]}=K[\{x^at^b\vert\, (a,b)\in
\mathbb{Z}^{n+1}\cap \mathbb{R}_+(I)\}].
\end{equation}
Hence, by Eqs.~(\ref{may6-06-2}) and (\ref{may6-06}), we get that
$R[It]$ is a normal domain if and 
only if  the following equality holds: 
$$
\mathbb{N}{\cal A}'=
\mathbb{Z}^{n+1}\cap\mathbb{R}_+(I).
$$
In geometric terms this means that $R[It]=\overline{R[It]}$ if
and only if ${\cal A}'$ is an integral Hilbert basis, that is, a 
Hilbert basis for the cone it generates. Rees algebras
and their integral closures are important 
objects of study in commutative algebra and geometry
\cite{bookthree}.   

A subset $C\subset X$ is a 
{\it minimal vertex cover\/} of the clutter $\cal C$ if: 
(i) every edge of $\cal C$ contains at least one vertex of $C$,  
and (ii) there is no proper subset of $C$ with the first 
property. If $C$ satisfies condition (i) only, then $C$ is 
called a {\it vertex cover\/} of $\cal C$. Let $C_1,\ldots,C_s$ be 
the minimal vertex covers of $\cal C$. The ideal of {\it vertex
covers\/}  
of ${\cal C}$ is the square-free monomial ideal
$$
I_c({\cal C})=(x^{u_1},\ldots,x^{u_s})\subset R,
$$
where $x^{u_k}=\prod_{x_i\in C_k}x_i$. The
clutter associated to $I_c({\cal C})$ is 
the {\it blocker\/} of $\cal C$, see \cite{cornu-book}. Notice that
the edges of the blocker are the minimal vertex covers of
$\cal C$.  

We now describe the content of the paper. A characterization
of perfect graphs--in terms of Rees cones--is
given (Proposition~\ref{myperfect-char}). We are able to
prove that $R[I_c({G})t]$ is normal if $G$ is a perfect graph
(Theorem~\ref{perfect-normal}) and that $R[I_c({G})t]$ is Gorenstein 
if $G$ is a perfect and unmixed graph (Corollary~\ref{perfect-gor}).
To show the normality of $R[I_c({G})t]$, we study when the system
$x\geq 0$; 
$xA\leq\mathbf{1}$ is TDI (Proposition~\ref{tdichar}), where TDI 
stands for Totally Dual Integral (see Section~\ref{perfectgraphs}). If
this system is TDI and the monomials in $F$ have the same degree, 
it is shown that $K[Ft]$ is an Ehrhart ring
(Proposition~\ref{jun28-06}). This is one of the results that will be
used 
in the proof of Theorem~\ref{perfect-normal}.

If $A$ is a balanced matrix, i.e., $A$ has no square submatrix of odd
order with exactly 
two $1$'s in 
each row and column, and $J=I_c({\cal C})$, then 
$R[It]=R_s(I)$ and 
$R[Jt]=R_s(J)$, see \cite{reesclu}. We complement these results
by showing that the Rees algebra of the dual $I^*$ of $I$ is normal
if $A$ is 
balanced (Proposition~\ref{dual-balanced-normal}).

By a result of Lyubeznik \cite{Lyu3}, $R_s(I({\cal C}))$ is a
$K$-algebra of finite  
type. Let $G$ be a graph. It is known that $R_s(I_c(G))$ is generated
as a $K$-algebra 
by monomials whose degree in $t$ is at most two
\cite[Theorem~5.1]{cover-algebras}, and one may even give an explicit
graph theoretical description of its minimal generators. Thus
$R_s(I_c(G))$ is well understood for graphs. In contrast, the minimal
set of generators of  $R_s(I(G))$ is very hard to describe in terms
of $G$ (see \cite{bahiano}). If $G$ is a perfect graph we compute the
integral Hilbert basis $\cal H$ of the Simis cone of $I(G)$ 
(see Definition~\ref{simis-cone-def} and
Theorem~\ref{symbolic-perfect}). 
Then, using 
that $R_s(I(G))$ is the semigroup ring of $\mathbb{N}{\cal H}$ over
$K$, we are able to prove that
$R_s(I(G))$ is generated as a $K$-algebra by monomials associated to
cliques of $G$ (Corollary~\ref{symbolic-perfect-coro}). 

Along the paper we introduce most of the 
notions that are relevant for our purposes. For unexplained
terminology and 
notation we refer to  
\cite{diestel,korte} and \cite{BHer,bookthree}. See \cite{cornu-book}
for additional information about clutters and perfect graphs.

\section{Perfect graphs, cones, and Rees algebras}\label{perfectgraphs}

We continue to use the notation and definitions used in the
introduction. Let $\mathfrak{p}_1,\ldots,\mathfrak{p}_s$ be the
minimal primes   
of $I({\cal C})$ and let 
$C_k=\{x_i\vert\, x_i\in\mathfrak{p}_k\}$ be the minimal vertex
cover of $\cal C$ 
that corresponds to $\mathfrak{p}_k$, see
\cite[Proposition~6.1.16]{monalg}.  
There is a
unique irreducible representation
$$
{\mathbb R}_+(I)=H_{e_1}^+\cap H_{e_2}^+\cap\cdots\cap
H_{e_{n+1}}^+\cap H_{\ell_1}^+\cap H_{\ell_2}^+\cap\cdots\cap
H_{\ell_r}^+
$$
such that each $\ell_k$ is in $\mathbb{Z}^{n+1}$, the non-zero
entries of 
each $\ell_k$ are relatively prime, and none of the closed 
halfspaces $H_{e_1}^+,\ldots,
H_{e_{n+1}}^+,H_{\ell_1}^+,\ldots,H_{\ell_r}^+$ can be
omitted from 
the intersection. Here $H_{a}^+$ denotes the closed halfspace 
$H_a^+=\{x\vert\, \langle
x,a\rangle\geq 0\}$ and $H_a$ stands for the hyperplane through the
origin with normal 
vector $a$, where $\langle\ ,\, \rangle$ denotes the standard 
inner product. According to \cite[Lemma~3.1]{normali} we
may always assume that $\ell_k=-e_{n+1}+\sum_{x_i\in C_k}e_i$ for
$1\leq k\leq s$. We shall be interested in the irreducible 
representation of the Rees cone of the ideal of vertex covers of a
perfect graph $G$ (see for instance
Proposition~\ref{myperfect-char}). 

Let $G$ be a simple graph with vertex set $X=\{x_1,\ldots,x_n\}$. In
what follows we shall always assume that $G$ has no isolated 
vertices. A {\it colouring\/} of the vertices of $G$ is an assignment
of colours to the vertices of $G$ in such a way that adjacent vertices
have distinct colours. The {\it chromatic number\/} of $G$ 
is the minimal number of colours in a colouring of $G$.
A graph is {\it perfect\/} if for every induced subgraph $H$, the
chromatic 
number of $H$ equals the size of the largest complete subgraph 
of $H$. The {\it complement\/} of $G$ is
denoted by $G'$. Recall that two vertices are adjacent in the graph
$G$ if and only if they are not adjacent in the graph $G'$.

Let $S$ be a subset of the vertices of $G$. The set $S$ is called 
{\it independent\/} if no two vertices of $S$ are adjacent. Notice
the following duality: $S$ is a maximal 
independent set of $G$ (with respect to inclusion) if and only if
$X\setminus S$ is a minimal vertex cover of $G$. We denote a complete
subgraph of $G$ with $r$ vertices  
by ${\cal K}_r$. The empty
set is regarded as an independent set whose incidence vector is the
zero vector. 

\begin{Theorem}{\rm (\cite[Theorem~16.14]{korte})}\label{chvatalperfect}
The following statements are equivalent\/{\rm :}
\begin{description}
\item{\rm (a)} $G$ is a perfect graph.\vspace{-1mm}
\item{\rm (b)} The complement of $G$ is perfect.\vspace{-1mm}
\item{\rm (c)} The independence polytope of $G$, i.e., the convex hull 
of the incidence vectors of the independent sets of $G$, is given
by\/{\rm :}
$$
\left\{(a_i)\in\mathbb{R}_+^{n}\vert\,
\textstyle\sum_{x_i\in\mathcal{K}_r}a_i\leq 
1;\ \forall\,  {\mathcal K}_r\subset G\right\}.
$$
\end{description}
\end{Theorem}
Below we express the perfection of $G$ in terms of a Rees cone. The 
next result is just a dual reinterpretation of part (c) above, which is
adequate to examine the normality and Gorensteiness 
of Rees algebras. We regard ${\cal K}_0$ as the empty set with zero
elements. A sum 
over an empty set is
defined to be $0$.

\begin{Proposition}\label{myperfect-char} Let $J=I_c(G)$ be the ideal
of vertex covers of $G$.  
Then $G$ is perfect if and only if the following equality holds
\begin{equation}\label{mundial06}
\mathbb{R}_+(J)=\left\{(a_i)\in\mathbb{R}^{n+1}\vert\,
\textstyle\sum_{x_i\in\mathcal{K}_r}a_i
\geq (r-1)a_{n+1};\ \forall\,  {\mathcal K}_r\subset
G\right\}.
\end{equation}
Moreover this is the irreducible representation of 
$\mathbb{R}_+(J)$ if $G$ is perfect. 
\end{Proposition}

\demo $\Rightarrow)$ The left hand side is contained in the right hand
side because  
any minimal vertex cover of $G$ contains at least $r-1$ vertices 
of any $\mathcal{K}_r$. For the reverse inclusion take a 
vector $a=(a_i)$ satisfying $b=a_{n+1}\neq 0$ and 
$$
\textstyle\sum_{x_i\in\mathcal{K}_r}a_i
\geq (r-1)b;\ \forall\,  {\mathcal K}_r\subset G
\ \Longrightarrow\ \textstyle\sum_{x_i\in\mathcal{K}_r}(a_i/b)
\geq r-1;\ \forall\,  {\mathcal K}_r\subset G.
$$
This implication follows because by making $r=0$ we get $b>0$. 
We may assume that $a_i\leq b$ for all $i$. Indeed 
if $a_i>b$ for some $i$, say $i=1$, then we can write
$a=e_1+(a-e_1)$. From the inequality
$$
\sum_{\substack{x_i\in{\cal K}_r\\ x_1\in{\cal K}_r}}a_i=a_1+
\sum_{\substack{x_i\in{\cal K}_{r-1}}}a_i\geq a_1+(r-2)b\geq 1+(r-1)b
$$
it is seen that $a-e_1$ belongs to the right hand side of
Eq.~(\ref{mundial06}). Thus, if necessary, we may apply this
observation again to $a-e_1$ and so on till we get that $a_i\leq b$ 
for all $i$. Hence, by  Theorem~\ref{chvatalperfect}(c), the vector 
$\gamma=\mathbf{1}-(a_1/b,\ldots,a_n/b)$
belongs to the independence polytope of $G$. Thus we can write
$$
\gamma=\lambda_1w_1+\cdots+\lambda_sw_s; \ \ \ (\lambda_i\geq 0;\,
\textstyle\sum_i\lambda_i=1),
$$
where $w_1,\ldots,w_s$ are incidence vectors of independent sets of
$G$. Hence
$$
\gamma=\lambda_1(\mathbf{1}-u_1')+\cdots+\lambda_s(\mathbf{1}-u_s'),
$$
where $u_1',\ldots,u_s'$ are incidence vectors of vertex covers of
$G$. Since any vertex cover contains a minimal one, for each $i$ we
can write $u_i'=u_i+\epsilon_i$, where $u_i$ is the incidence
vector of a minimal vertex cover of $G$ and $\epsilon_i\in\{0,1\}^n$.
Therefore   
\begin{eqnarray*}
\lefteqn{\mathbf{1}-\gamma=\lambda_1u_1'+\cdots+\lambda_su_s'
\Longrightarrow}\ \ \ \ \ \ \ \ \ \ \ \ \ \ \ \ \ \ \ \ \ \ \ \ \ \ \ \ \
\ \ \ \ \ \ \ \ \ \ \ \ \ \ \ \ \ \ \ \ \ \ \ \ \ \ \ 
\ \ \ \ \ \ \ \ \ \ \ \ \ \ \ \ \ \ \ \ \ \ \ \ \\  
a=b\lambda_1(u_1,1)+\cdots+b\lambda_s(u_s,1)+b\lambda_1\epsilon_1+\cdots+
b\lambda_s\epsilon_s,
\end{eqnarray*}
Thus $a\in\mathbb{R}_+(J)$. 
If $b=0$, clearly $a\in\mathbb{R}_+(J)$. Hence we get 
equality in Eq.~(\ref{mundial06}), as required. The converse follows using
similar arguments. 

To finish the proof it suffices to show that the set
$$
F=\{(a_i)\in\mathbb{R}^{n+1}\vert\,
\textstyle\sum_{x_i\in\mathcal{K}_r}a_i=(r-1)a_{n+1}\}\cap\mathbb{R}_+(J)
$$
is a facet of $\mathbb{R}_+(J)$. If $\mathcal{K}_r=\emptyset$, then
$r=0$ and $F=H_{e_{n+1}}\cap\mathbb{R}_+(J)$, which is clearly a
facet because  
$e_1,\ldots,e_n\in F$. If $r=1$, then $F=H_{e_i}\cap\mathbb{R}_+(J)$ 
for some $1\leq i\leq n$, which is a facet because $e_j\in F$ for
$j\notin\{i,n+1\}$ and there is at least 
one minimal vertex cover of $G$ not containing $x_i$. We may assume that
$X'=\{x_1,\ldots,x_r\}$ is the vertex set of
${\cal K}_r$ and $r\geq 2$. For each
$1\leq i\leq r$ there is a minimal vertex cover $C_i$ of $G$ not 
containing $x_i$. Notice that $C_i$ contains $X'\setminus\{x_i\}$.
Let $u_i$ be the incidence vector 
of $C_i$. Since the rank of $u_1,\ldots,u_r$ is $r$, it follows that
the set 
$$
\{(u_1,1),\ldots,(u_r,1),e_{r+1},\ldots,e_n\}
$$
is a linearly independent set contained in $F$, i.e., $\dim(F)=n$.
Hence $F$ is a facet of $\mathbb{R}_+(J)$ because the hyperplane that
defines $F$ is a supporting hyperplane. \QED

\medskip

There are computer programs that determine 
the irreducible representation of a Rees cone \cite{B}. 
Thus we may use Proposition~\ref{myperfect-char} to 
determine whether a given graph is perfect, and in 
the process we may also determine its complete subgraphs. However this
proposition is useful mainly for theoretical reasons. A direct
consequence of this 
result (Lemma~\ref{mundial09-jun-06}(b) below) will be used to
prove one of our main results (Theorem~\ref{perfect-normal}).

Let $S$ be a set of vertices of a graph $G$, the {\it induced
subgraph\/} $\langle S\rangle$ is the maximal subgraph of $G$ with
vertex set $S$. A {\it clique\/} of a graph $G$ is a subset of the set
of vertices that induces 
a complete subgraph. We will also call a complete subgraph of $G$ 
a clique. The {\it support\/} of 
$x^a=x_1^{a_1}\cdots x_n^{a_n}$ is ${\rm supp}(x^a)= \{x_i\, |\,
a_i>0\}$. If $a_i\in\{0,1\}$ for all $i$, $x^a$ is called a {\it
square-free\/} monomial. We regard the empty set as an
independent set with zero elements.

\begin{Lemma}\label{mundial09-jun-06} {\rm (a)}
$I_c(G')=(\{x^a\vert\, X\setminus{\rm 
supp}(x^a)\mbox{ is a maximal clique of }G\})$.

\item{\rm (b)} If $G$ is perfect and $J'=I_c(G')$, then $\mathbb{R}_+(J')$ 
is equal to 
$$\left\{(a_i)\in\mathbb{R}^{n+1}\vert\,
\textstyle\sum_{x_i\in S}a_i
\geq (|S|-1)a_{n+1};\ \forall\, S\mbox{ independent set of }G \right\}.
$$
\end{Lemma}

\demo (a) Let $x^a\in R$ and let $S={\rm supp}(x^a)$. Then $x^a$ is a
minimal generator of $I_c(G')$ if and only if $S$ is a minimal vertex cover
of $G'$ if and only if $X\setminus S$ is a maximal independent set of
$G'$ if and only if $\langle X\setminus S\rangle$ is a maximal complete 
subgraph of $G$. Thus the equality holds. (b) By
Theorem~\ref{chvatalperfect} the graph $G'$ is perfect. Hence 
the equality follows from Proposition~\ref{myperfect-char}. \QED 

\medskip

Let $A$ be an integral matrix. 
The system $x\geq 0;\ xA\leq\mathbf{1}$ is called 
{\it totally dual integral\/} (TDI) if the minimum in 
the LP-duality equation
\begin{equation}\label{jun6-2-03}
{\rm max}\{\langle \alpha,x\rangle \vert\, x\geq 0; xA\leq \mathbf{1}\}=
{\rm min}\{\langle y,\mathbf{1}\rangle \vert\, y\geq 0; Ay\geq\alpha\} 
\end{equation}
has an integral optimum solution $y$ for each integral vector $\alpha$ with 
finite minimum. 

An incidence matrix $A$ of a clutter is called {\it
perfect\/} if the polytope defined by the system $x\geq 0;\
xA\leq\mathbf{1}$ is integral, i.e., it has only integral vertices. 
The {\it vertex-clique matrix\/} of a graph $G$ is the $\{0,1\}$-matrix
whose rows are indexed by the vertices of $G$ and whose columns are
the incidence vectors of the maximal cliques of $G$.

\begin{Theorem}[\rm\cite{lovasz},\cite{chvatal}]\label{lovasz-chvatal}
Let $A$ be the 
incidence matrix of a clutter. Then the
following are equivalent{\rm :}
\begin{description}
\item{\rm (a)} The system $x\geq 0;\ xA\leq\mathbf{1}$ is {\rm
TDI}.\vspace{-1mm}
\item{\rm (b)} $A$ is perfect.\vspace{-1mm}
\item{\rm (c)} $A$ is the vertex-clique matrix of a perfect graph.  
\end{description}
\end{Theorem}

\begin{Proposition}\label{tdichar} 
Let $A$ be an $n\times q$ matrix
with entries in $\mathbb{N}$ and let $v_1,\ldots,v_q$ be its column
vectors. Then the system 
$x\geq 0;\ xA\leq\mathbf{1}$ is {\rm TDI} if and only if 
\begin{description}
\item{\rm (i)\ } the polyhedron $\{x\vert\, x\geq 0;\,
xA\leq\mathbf{1}\}$ is integral, and \vspace{-1mm}
\item{\rm (ii)} $\mathbb{R}_+{\cal
B}\cap\mathbb{Z}^{n+1}=\mathbb{N}{\cal B}$, where ${\cal
B}=\{(v_1,1),\ldots,(v_q,1),-e_1,\ldots,-e_n\}$. 
\end{description}
\end{Proposition}
\demo $\Rightarrow$) By \cite[Corollary~22.1c]{Schr} we get that (i)
holds. To prove (ii) take $(\alpha,b)\in\mathbb{R}_+{\cal
B}\cap\mathbb{Z}^{n+1}$, where $\alpha\in\mathbb{Z}^n$ and 
$b\in\mathbb{Z}$. By hypothesis the minimum in Eq.~(\ref{jun6-2-03})
has an integral optimum solution $y=(y_i)$ such that
$|y|=\langle{y},\mathbf{1}\rangle\leq b$.  
Since $y\geq 0$ and $\alpha\leq Ay$ we can write
$$
\begin{array}{l} 
\ \ \ \ \ \alpha=y_1v_1+\cdots+y_qv_q-\delta_1e_1-\cdots-\delta_ne_n\ \ \ 
(\delta_i\in\mathbb{N})\ \ \Longrightarrow\\ 
(\alpha,b)=y_1(v_1,1)+\cdots+y_{q-1}(v_{q-1},1)+(y_q+b-|y|)(v_q,1)
-(b-|y|)v_q-\delta,
\end{array}
$$
where $\delta=(\delta_i)$. As the entries of $A$ are in $\mathbb{N}$,
the vector $-v_q$ can be written as a non-negative integer
combination of $-e_1,\ldots,-e_n$. Thus $(\alpha,b)\in\mathbb{N}{\cal B}$.
This proves (ii).

$\Leftarrow$) Assume that the system $x\geq 0$; $xA\leq\mathbf{1}$ is
not TDI.  
Then there exists an $\alpha_0\in {\mathbb Z}^n$ such that if ${y}_0$ 
is an optimal solution of the linear program:
\begin{equation}\label{brasil-ghana}
\min\{\langle y,\mathbf{1}\rangle\vert\;y\geq{0};\; Ay\geq{\alpha}_0\},
\end{equation}
then ${y}_0$ is not integral. We claim that also the optimal value 
$|y_0|=\langle{y}_0,\mathbf{1}\rangle$ of this linear program is not
integral. If  
$|y_0|$ is integral, 
then $(\alpha_0,|y_0|)$ is in $\mathbb{Z}^{n+1}\cap\mathbb{R}_+{\cal
B}$. Hence by (ii), we get that $(\alpha_0,|y_0|)$ is in 
$\mathbb{N}{\cal B}$, but this readily yields that the linear program
of Eq.~(\ref{brasil-ghana}) has an integral optimal solution, a
contradiction. This completes the  
proof of the claim. Consider the dual linear program:
$$
\max\{\langle x,\alpha_0\rangle\vert \; x \geq{0},\;  xA\leq \mathbf{1}\}.
$$
Its optimal value is attained at a 
vertex $x_0$ of $\{x\vert\, x\geq 0;\,
xA\leq\mathbf{1}\}$. Then by LP duality we 
get $\langle x_0,\alpha_0\rangle=|{y}_0|\notin {\mathbb Z}$. Hence $x_0$ is not 
integral, a contradiction to the integrality of $\{x\vert\, x\geq 0;\,
xA\leq\mathbf{1}\}$. \QED

\begin{Remark}\rm If $A$ is a matrix with entries in $\mathbb{Z}$
satisfying (i) and (ii), then the system $x\geq 0;\ xA\leq\mathbf{1}$
is {\rm TDI}.
\end{Remark}

Let $v_1,\ldots,v_q$ be a set of points in $\mathbb{N}^n$ and let 
$P={\rm conv}(v_1,\ldots,v_q)$. The
{\it Ehrhart ring\/} of the lattice polytope $P$ is the $K$-subring
of $R[t]$  given by
$$
A(P)=K[\{x^at^b\vert\, a\in bP\cap\mathbb{Z}^n\}]. 
$$

\begin{Proposition}\label{jun28-06} Let $A$ be a perfect matrix
with column vectors $v_1,\ldots,v_q$. If there is
$x_0\in\mathbb{R}^n$ such that 
all the entries of $x_0$ are positive and $\langle v_i,x_0\rangle=1$ 
for all $i$, then 
$A(P)=K[x^{v_1}t,\ldots,x^{v_q}t]$.
\end{Proposition}

\demo Let $x^at^b\in A(P)$. Then we can write 
$(a,b)=\sum_{i=1}^q\lambda_i(v_i,1)$, 
where $\lambda_i\geq 0$ for all $i$. 
Hence $\langle a,x_0\rangle=b$. By Theorem~\ref{lovasz-chvatal} the
system $x\geq 0$; 
$xA\leq\mathbf{1}$ is TDI. Hence applying 
Proposition~\ref{tdichar}(ii) we have:
$$ 
(a,b)=\eta_1(v_1,1)+\cdots+\eta_q(v_q,1)-\delta_1e_1-\cdots-\delta_ne_n
\ \ \  (\eta_i\in\mathbb{N};\, \delta_i\in\mathbb{N}).
$$
Consequently $b=\langle
a,x_0\rangle=b-\delta_1\langle x_0,e_1\rangle-\cdots-\delta_n\langle
x_0,e_n\rangle$. Using that $\langle x_0,e_i\rangle>0$ for all $i$,  
we conclude that $\delta_i=0$ for all $i$, 
i.e., $x^at^b\in K[x^{v_1}t,\ldots,x^{v_q}t]$.  \QED  

\medskip

Recall that the clutter $\cal C$ (or the edge ideal
$I({\cal C})$) is called {\it unmixed\/} if all the minimal vertex
covers of $\cal C$  
have the same cardinality. 

\begin{Corollary}
If $G$ is a perfect 
unmixed graph and $v_1,\ldots,v_q$ are the
incidence vectors of the maximal independent sets of $G$, 
then $K[x^{v_1}t,\ldots,x^{v_q}t]$ is normal. 
\end{Corollary}

\demo The minimal vertex covers of $G$ are exactly 
the complements of the maximal independent sets of $G$. Thus 
$|v_i|=d$ for all $i$, where $d=\dim(R/I(G))$. On the other hand 
the maximal independent sets of $G$ are exactly the maximal
cliques of $G'$. Thus, by Theorem~\ref{lovasz-chvatal} and
Proposition~\ref{jun28-06}, the subring
$K[x^{v_1}t,\ldots,x^{v_q}t]$ is  
an Ehrhart ring, and consequently it is normal.  \QED  

\medskip

Let $\cal C$ be a clutter and let $A$ be its incidence matrix. The 
clutter $\cal C$ satisfies the {\it max-flow min-cut\/} (MFMC)
property if both sides 
of the LP-duality equation
$$
{\rm min}\{\langle \alpha,x\rangle \vert\, x\geq 0; xA\geq \mathbf{1}\}=
{\rm max}\{\langle y,\mathbf{1}\rangle \vert\, y\geq 0; Ay\leq\alpha\} 
$$
have integral optimum solutions $x$ and $y$ for each non-negative
integral vector $\alpha$, see \cite{cornu-book}. Let $I$ be the edge
ideal of $\cal C$. Closely related 
to ${\mathbb R}_+(I)$ is the {\it set covering polyhedron\/}: 
\[
Q(A)=\{x \in \mathbb{R}^n  \, | \; x \geq{0}, \;
xA\geq{\mathbf 1} \},
\]
see \cite[Theorem~3.1]{reesclu}. Its integral
vertices 
are precisely the incidence vectors of the minimal vertex covers of
$\cal C$ \cite[Proposition~2.2]{reesclu}. 

\begin{Corollary}
Let $\cal C$ be a clutter and
let $A$ be its incidence 
matrix. If all the edges of $\cal C$ have the same cardinality and 
the polyhedra 
$$
\{x\vert\, x\geq 0;\,
xA\leq\mathbf{1}\}\ \mbox{ and }\ \{x\vert\, x\geq 0;\,
xA\geq\mathbf{1}\}
$$
are integral, then $\cal C$ has the max-flow min-cut property. 
\end{Corollary}

\demo By \cite[Proposition~4.4 and Theorem~4.6]{reesclu} we have that
$\cal C$ has the max-flow min-cut property if and only if $Q(A)$ is integral 
and $K[x^{v_1}t,\ldots,x^{v_q}t]=A(P)$, where $v_1,\ldots,v_q$ are the
column vectors of $A$ and $P={\rm conv}(v_1,\ldots,v_q)$. Thus 
the result follows from Proposition~\ref{jun28-06}. \QED

\medskip

The {\it clique clutter\/} of a graph $G$, denoted
by ${\rm cl}(G)$, is the clutter on $V(G)$ whose edges are the 
maximal cliques of $G$. 

\begin{Theorem}\label{perfect-normal} If $G$ is a perfect graph, then
$R[I_c(G)t]$ is normal. 
\end{Theorem}

\demo Let $G'$ be the complement of $G$ and let $J'=I_c(G')$. Since
$G'$ is perfect it suffices to prove that $R[J't]$ is normal. 

Case (A): Assume that all the maximal cliques of $G$ have the same
number of elements. Let $F=\{x^{v_1},\ldots,x^{v_q}\}$ be the set of
monomials of $R$ 
whose support is a maximal clique of $G$. We set 
$F'=\{x^{w_1},\ldots,x^{w_q}\}$, where $x^{w_i}=x_1\cdots
x_n/x^{v_i}$. By Lemma~\ref{mundial09-jun-06}(a) we have 
$J'=(F')$. Consider the matrices 
$$
B=\left(\hspace{-1mm}
\begin{array}{ccc}
v_1&\cdots&v_q\\ 
1&\cdots &1
\end{array}\hspace{-1mm}
\right)\ \mbox{ and }\ B'=\left(\hspace{-1mm}
\begin{array}{ccc}
w_1&\cdots&w_q\\ 
1&\cdots &1
\end{array}\hspace{-1mm}
\right),
$$
where the $v_i$'s and $w_j$'s are regarded as column vectors. Using 
the last row of $B$ as a pivot it is seen that $B$ is equivalent over
$\mathbb{Z}$ to $B'$. Let $A$ be the incidence matrix of ${\rm
cl}(G)$, the clique clutter of $G$, whose
columns are $v_1,\ldots,v_q$. As the matrix $A$ is perfect, by 
Proposition~\ref{jun28-06}, we obtain that $K[Ft]=A(P)$, where 
$A(P)$ is the Ehrhart ring of $P={\rm conv}(v_1,\ldots,v_q)$. In
particular $K[Ft]$ is normal because Ehrhart rings are normal. 
According to \cite[Theorem~3.9]{ehrhart} we have that $K[Ft]=A(P)$  
if and only if $K[Ft]$ is normal and $B$ diagonalizes over $\mathbb{Z}$  
to an ``identity'' matrix. Consequently the matrix $B'$ diagonalizes to an
identity matrix along with $B$.  Since 
the rings $K[F't]$ and $K[Ft]$ 
are isomorphic, we get that $K[F't]$ is normal. Thus, again by 
\cite[Theorem~3.9]{ehrhart}, we obtain the equality $K[F't]=A(P')$, where
$A(P')$ is the Ehrhart ring of $P'={\rm conv}(w_1,\ldots,w_q)$. 
Let $H_a^+$ be any of the halfspaces that occur 
in the irreducible representation of the Rees cone $\mathbb{R}_+(J')$. By 
Lemma~\ref{mundial09-jun-06}(b) the first $n$ entries of $a$ are
either $0$ or $1$. Hence by 
\cite[Proposition~4.2]{reesclu} we get the equality 
$$
A(P')[x_1,\ldots,x_n]=\overline{R[J't]}. 
$$
Therefore $R[J't]=K[F't][x_1,\ldots,x_n]=
A(P')[x_1,\ldots,x_n]=\overline{R[J't]}$, that is, 
$R[J't]$ is normal. 

Case (B): Assume that not all the maximal cliques of $G$ have the 
same number of elements. Let $C$ be a maximal clique of $G$ of lowest
size and let $w$ be its incidence vector. For simplicity of notation 
assume that $C=\{x_1,\ldots,x_r\}$. Let $z=x_{n+1}\notin V(G)$ be a
new vertex.  We construct a new graph $H$ as follows. Its vertex set 
is $V(H)=V(G)\cup\{z\}$ and its edge set is 
$$
E(H)=E(G)\cup\{\{z,x_1\},\ldots,\{z,x_r\}\}.
$$
Notice that $C\cup\{z\}$ is the only maximal clique of $H$ containing
$z$. 
Thus it is seen that the edges of the clique clutter of $H$ are related to
those of the clique clutter of $G$ as follows:
$$
E({\rm cl}(H))=(E({\rm cl}(G))\setminus\{C\})\cup\{C\cup\{z\}\}.
$$
From the proof of \cite[Proposition~5.5.2]{diestel} it follows that if
we paste together $G$ and the complete subgraph induced by
$C\cup\{z\}$ 
along the complete subgraph induced by $C$ we
obtain a perfect graph, i.e., $H$ is perfect. This construction is
different from the famous Lov\'asz replication of a vertex, 
as explained in \cite[Lemma~3.3]{cornu-book}. The contraction of 
${\rm cl}(H)$ at $z$, denoted by ${\rm cl}(H)/z$, is the clutter of
minimal elements 
of $\{S\setminus\{z\}\vert\, S\in {\rm cl}(H)\}$. In our case we 
have ${\rm cl}(H)/z={\rm cl}(G)$, i.e., ${\rm cl}(G)$ is a minor 
of ${\rm cl}(H)$ obtained by contraction. By successively adding 
new vertices $z_1=z,z_2,\ldots,z_r$, following the construction above, we
obtain a perfect  
graph $H$ whose maximal cliques have the same size and such that 
${\rm cl}(G)$ is a minor of ${\rm cl}(H)$ obtained by contraction 
of the vertices $z_1,\ldots,z_s$. By case (A) we obtain that the ideal
$L=I_c(H')$ of minimal vertex covers of $H'$ is normal. Since $L$ is 
generated by all the square-free monomials $m$ of $R[z_1,\ldots,z_s]$
such that 
$V(H)\setminus {\rm supp}(m)$ is a maximal clique of $H$, it follows 
that $J'$ is obtained from $L$ by making $z_i=1$ for all $i$.
Hence $R[J't]$ is normal because the normality property of 
Rees algebras of edge ideals is closed under taking 
minors \cite[Proposition~4.3]{normali}. \QED 

\begin{Example}\rm If $G$ is a pentagon, then the Rees algebra of
$I_c(G)$ is normal and 
$G$ is not perfect. 
\end{Example}

\begin{Corollary}\label{perfect-gor} If $G$ is perfect and unmixed, 
then $R[I_c(G)t]$ is a Gorenstein standard graded $K$-algebra.
\end{Corollary}

\demo Let $g$ be the height of the edge ideal $I(G)$ and let
$J=I_c(G)$. By assigning $\deg(x_i)=1$ and    
$\deg(t)=-(g-1)$, the Rees algebra $R[Jt]$ becomes a graded
$K$-algebra generated by monomials of  degree $1$. The Rees ring
$R[Jt]$ is a normal     
domain by Theorem~\ref{perfect-normal}. Then according to a formula
of  Danilov-Stanley \cite[Theorem~6.3.5]{BHer}  its canonical module
is the ideal     
of $R[Jt]$ given by   
$$   
\omega_{R[Jt]}=(\{x_1^{a_1}\cdots x_n^{a_n}   
t^{a_{n+1}}\vert\, a=(a_i)   
\in\mathbb{R}_+(J)^{\rm o}\cap\mathbb{Z}^{n+1}\}),   
$$   
where $\mathbb{R}_+(J)^{\rm o}$ denotes the topological    
interior of the Rees cone of $J$. By a result of Hochster \cite{Ho1} 
the ring $R[Jt]$ is Cohen-Macaulay. Using Eq.~(\ref{mundial06}) it is
seen that the vector     
$(1,\ldots,1)$ is in the interior of the Rees cone, i.e., $x_1\cdots
x_nt$ belongs to $\omega_{R[Jt]}$. Take an    
arbitrary monomial    
$x^at^b=x_1^{a_1}\cdots x_n^{a_n}t^b$ in the ideal $\omega_{R[Jt]}$,  
that is $(a,b)\in\mathbb{R}_+(J)^{\rm o}$. Hence the vector $(a,b)$
has positive integer entries   
and    
satisfies   
\begin{eqnarray}   
& \sum_{x_i\in \mathcal{K}_r}a_i\geq(r-1)b+1&\label{nov9-03}  
\end{eqnarray}   
for every complete subgraph $\mathcal{K}_r$ of $G$. If $b=1$, clearly
$x^at^b$ is a multiple of $x_1\cdots x_nt$. 
Now assume $b\geq 2$. Using the normality of    
$R[Jt]$ and Eqs.~(\ref{mundial06})    
and (\ref{nov9-03})   
it follows that the monomial $m=x_1^{a_1-1}\cdots x_n^{a_n-1}t^{b-1}$    
belongs to $R[Jt]$. Since $x^at^b=mx_1\cdots x_nt$, we obtain that
$\omega_{R[Jt]}$ is     
generated by $x_1\cdots x_nt$ and thus $R[Jt]$ is a    
Gorenstein ring.\QED

\medskip

A graph ${G}$ is {\it chordal\/} if every
cycle of ${G}$ of length $n\geq 4$ has a chord. A {\it chord} of a
cycle is an edge joining two non adjacent vertices of the cycle.

\begin{Corollary}
If $J$ is a Cohen-Macaulay square-free monomial ideal of height two, 
then $R[Jt]$ is normal.  
\end{Corollary}

\demo Consider the graph $G$ whose edges are the pairs $\{x_i,x_j\}$
such that $(x_i,x_j)$ is a minimal prime of $J$. Notice that
$J=I_c(G)$. By \cite[Theorem~6.7.13]{monalg}, the ideal $I_c(G)$ is
Cohen-Macaulay if and only if $G'$ is a chordal graph. Since chordal
graphs are perfect 
\cite[Proposition~5.5.2]{diestel}, we obtain that $G'$ is perfect. 
Thus $G$ is a perfect graph by
Theorem~\ref{chvatalperfect}. Applying Theorem~\ref{perfect-normal}
we conclude that $R[Jt]$ is normal. \QED

\medskip

Recall that a matrix with $\{0,1\}$-entries is called {\it balanced\/} if
$A$ has no square submatrix of odd order with exactly
two $1$'s in each row and column, 

\begin{Proposition}\label{dual-balanced-normal} 
Let $A$ be a $\{0,1\}$-matrix with column
vectors $v_1,\ldots,v_q$ and let $w_i=\mathbf{1}-v_i$. If $A$ is
balanced, 
then the Rees algebra of 
$I^*=(x^{w_1},\ldots,x^{w_q})$ is a normal domain.  
\end{Proposition}

\demo According to \cite{berge-balanced}, \cite[Corollary~83.1a(vii),
p.1441]{Schr2} $A$ is balanced if and only if every
submatrix of $A$ is perfect. By adjoining rows of unit vectors to $A$
and since the normality property of edge ideals is closed under taking
minors \cite[Proposition~4.3]{normali} we may assume that $|v_i|=d$
for all 
$i$. By Theorem~\ref{lovasz-chvatal} there is a perfect graph $G$ such
that $A$ is the vertex-clique matrix of $G$. Thus following the first
part of the proof of 
Theorem~\ref{perfect-normal}, we obtain that $R[I^*t]$ is normal. \QED

\medskip

Consider the ideals $I=(x^{v_1},\ldots,x^{v_q})$ and 
$I^*=(x^{w_1},\ldots,x^{w_q})$. Following the terminology of matroid
theory we call $I^*$ the {\it
dual\/} of $I$. Notice the following duality. If $A$ is the
vertex-clique matrix of a graph $G$,  
then $I^*$ is precisely the ideal of vertex covers of $G'$.  

\section{Symbolic Rees algebras of  edge ideals}

Let $G$ be a graph with vertex set $X=\{x_1,\ldots,x_n\}$ and let 
$I=I(G)$ be its edge ideal \cite[Chapter~6]{monalg}. 
The main purpose of this section is to study the symbolic Rees 
algebra of $I$ and the Simis cone of $I$ when $G$ is a perfect
graph. We show that the 
cliques of a perfect graph $G$ completely determine both the Hilbert
basis 
of the Simis cone and the symbolic Rees algebra of $I(G)$.  

\begin{Definition}\label{simis-cone-def}\rm 
Let $C_1,\ldots,C_s$ be the 
minimal vertex covers of $G$. The {\it symbolic Rees cone\/} 
or {\it Simis cone} of $I$ is the
rational polyhedral cone: 
$$
{\rm Cn}(I)=H_{e_1}^+\cap\cdots\cap H_{e_{n+1}}^+\cap
H_{(u_1,-1)}^+\cap\cdots\cap H_{(u_s,-1)}^+,
$$
where $u_k=\sum_{x_i\in C_k}e_i$ for $1\leq k\leq s$. 
\end{Definition}

Simis cones were
introduced in \cite{normali} to study symbolic Rees algebras of
square-free monomial ideals. If $\cal H$ is an integral Hilbert basis
of ${\rm Cn}(I)$, then $R_s(I(G))$ equals
$K[\mathbb{N}{\cal H}]$, the semigroup ring of $\mathbb{N}{\cal H}$
(see \cite{normali}). This result is interesting because it 
allows us to compute the 
minimal generators of $R_s(I(G))$ using Hilbert basis. Next we
describe $\cal H$ when $G$ is perfect. 

\begin{Theorem}\label{symbolic-perfect} Let 
$\omega_1,\ldots,\omega_p$ be the incidence vectors of the non-empty
cliques of a perfect graph $G$ and let 
$$
{\cal H}=\{(\omega_1,|\omega_1|-1),\ldots,(\omega_p,|\omega_p|-1)\}.
$$
Then $\mathbb{N}{\cal H}={\rm Cn}(I)\cap\mathbb{Z}^{n+1}$, where
$\mathbb{N}{\cal H}$ is the subsemigroup of $\mathbb{N}^{n+1}$ generated
by $\cal H$, that is, 
$\cal H$ is the integral Hilbert basis of ${\rm Cn}(I)$. 
\end{Theorem}

\demo The inclusion $\mathbb{N}{\cal H}\subset{\rm
Cn}(I)\cap\mathbb{Z}^{n+1}$ is clear because each clique of size $r$
intersects any minimal vertex cover in at least $r-1$ vertices. Let us
show the reverse inclusion. Let $(a,b)$ be a minimal generator 
of ${\rm Cn}(I)\cap\mathbb{Z}^{n+1}$, where $0\neq
a=(a_i)\in\mathbb{N}^n$ and 
$b\in\mathbb{N}$. Then
\begin{equation}\label{apr04-07}
\textstyle\sum_{x_i\in C_k}a_i=\langle a,u_k\rangle\geq b,
\end{equation}
for all $k$. If $b=0$ or $b=1$, then $(a,b)=e_i$ for some
$i\leq n$ or $(a,b)=(e_i+e_j,1)$ for some edge $\{x_i,x_j\}$
respectively. In both cases $(a,b)\in{\cal H}$. Thus we  
may assume that $b\geq 2$ and $a_j\geq 1$ for some $j$. Using 
Eq.~(\ref{apr04-07}) we obtain
\begin{equation}\label{aug24-06}
\sum_{x_i\in C_k}a_i+\hspace{-2mm}\sum_{x_i\in X\setminus C_k}
\hspace{-2mm}a_i=|a|\geq b+
\hspace{-2mm}\sum_{x_i\in X\setminus C_k}\hspace{-2mm}
a_i=b+\langle\mathbf{1}-u_k,a \rangle,
\end{equation}
for all $k$, where $X=\{x_1,\ldots,x_n\}$ is the vertex set of $G$.
Set $c=|a|-b$. 
Notice that $c\geq 1$ because $a\neq 0$. Indeed if $c=0$, from 
Eq.~(\ref{aug24-06}) we get $\sum_{\scriptstyle x_i\in X\setminus
C_k}a_i=0$ for all 
$k$, i.e., $a=0$, a contradiction. Consider the vertex-clique matrix of
$G'$:
$$
A'=\left(\mathbf{1}-u_1\cdots \mathbf{1}-u_s
\right),
$$
where $\mathbf{1}-u_1,\ldots, \mathbf{1}-u_s$ are regarded as column
vectors. From Eq.~(\ref{aug24-06}) we get $(a/c)A'\leq\mathbf{1}$.
Hence by Theorem~\ref{chvatalperfect}(c) we obtain that $a/c$ belongs
to ${\rm conv}(\omega_0,\omega_1,\ldots,\omega_p)$, where
$\omega_0=0$, i.e., we can write
$a/c=\lambda_0\omega_0+\cdots+\lambda_p\omega_p$, where 
$\lambda_i\geq 0$ for all $i$ and $\sum_{i}\lambda_i=1$. Thus we can
write
$$
(a,c)=c\lambda_0(\omega_0,1)+\cdots+c\lambda_p(\omega_p,1).
$$
Using Theorem~\ref{lovasz-chvatal}(a) it is not hard to see that the 
subring $K[\{x^{\omega_i}t\vert 0\leq i\leq p\}]$ is normal. 
Hence there are $\eta_0,\ldots,\eta_p$ in $\mathbb{N}$ such that 
$$
(a,c)=\eta_0(\omega_0,1)+\cdots+\eta_p(\omega_p,1).
$$ 
Thus $|a|=\eta_0|\omega_0|+\cdots+\eta_p|\omega_p|$ and
$c=\eta_0+\cdots+\eta_p=|a|-b$, consequently:
$$
(a,b)=\eta_0(\omega_0,|\omega_0|-1)+\eta_1(\omega_1,|\omega_1|-1)
+\cdots+\eta_p(\omega_p,|\omega_p|-1).
$$ 
Notice that
there is $u_\ell$ such that $\langle 
a,u_\ell\rangle=b$; otherwise since $a_j\geq 1$, by
Eq.~(\ref{apr04-07}) the
vector $(a,b)-e_j$ would be in ${\rm Cn}(I)\cap\mathbb{Z}^{n+1}$, 
contradicting the minimality of $(a,b)$.
Therefore from the equality
$$
0=\langle(a,b),(u_\ell,-1)\rangle=\eta_0+\textstyle\sum_{i=1}^p
\eta_i\langle(\omega_i,|\omega_i|-1),(u_\ell,-1)\rangle
$$ 
we conclude that $\eta_0=0$, i.e., $(a,b)\in\mathbb{N}{\cal H}$, as
required. \QED
 
\begin{Corollary}\label{symbolic-perfect-coro} 
If $G$ is a perfect graph, then 
$$
R_s(I(G))=K[x^at^r\vert\, 
x^a\mbox{ is square-free };\, 
\langle{\rm supp}(x^a)\rangle={\cal K}_{r+1};\, 0\leq
r<n].
$$
\end{Corollary}

\demo Let $K[\mathbb{N}{\cal H}]$ be the semigroup ring with
coefficients in $K$ of the semigroup $\mathbb{N}{\cal H}$. By
\cite[Theorem~3.5]{normali} we have 
the equality $R_s(I(G))=K[\mathbb{N}{\cal H}]$, thus the formula
follows from  
Theorem~\ref{symbolic-perfect}. \QED

\begin{Corollary}[\rm \cite{bahiano}]
If $G$ is a complete graph, then 
$$
R_s(I(G))=K[x^at^r\vert\, 
x^a\mbox{ is square-free };\, \deg(x^a)={r+1};\, r\geq 0].
$$
\end{Corollary}

\bibliographystyle{plain}

\end{document}